\documentclass[10pt]{article}
\usepackage{latexsym}
\usepackage{amsfonts}
\usepackage{enumerate}
\usepackage{multicol}
\usepackage{graphicx}
\usepackage{amssymb}
\usepackage{amsmath}
\usepackage{epic}

\title{Comments on Y.~O.~Hamidoune's Paper ``Adding Distinct Congruence Classes''  \\[.4in]}

\author{B\'{e}la Bajnok \\[.1in] Department of Mathematics, Gettysburg College \\
Gettysburg, PA 17325-1486 USA \\E-mail:  bbajnok@gettysburg.edu \\[.4in]}

\date{June 18, 2015}

\newtheorem{thm}{Theorem}

\newtheorem{lem}[thm]{Lemma}

\newtheorem{prop}[thm]{Proposition}

\begin{document}

\maketitle

\begin{abstract}
The main result in Y.~O.~Hamidoune's paper ``Adding Distinct Congruence Classes'' ({\em Combin.~Probab.~Comput.}~{\bf 7} (1998) 81-87) is as follows:  
If $S$ is a generating subset of a cyclic group $G$ such that $0 \not \in S$ and $|S| \geq 5$, then the number of sums of the subsets of $S$ is at least $\min (|G|, 2|S| )$.  Unfortunately, argument of the author, who, sadly, passed away in 2011, relies on a lemma whose proof is incorrect; in fact, the lemma is false for all cyclic groups of even order.  In this short note we point out this mistake, correct the proof, and discuss why the main result is actually true for all finite abelian groups.

\end{abstract}


\thispagestyle{empty}

Let $G$ be a finite abelian group, written additively.  For a positive integer $h$ and a subset $A$ of $G$, we let $h \hat{\;} A$ denote the set of sums of the $h$-subsets of $A$: 
$$h \hat{\;} A = \{ \Sigma_{b \in B} b \mid B \subseteq A, |B|=h\}.$$   Additionally, we let $\Sigma A$ denote the set of all nonempty subset sums of $A$:
$$\Sigma A = \cup_{h=1}^{|A|} h \hat{\;} A = \{ \Sigma_{b \in B} b \mid \emptyset \neq B \subseteq A\}.$$  If $G$ is cyclic and of order $m$, we identify it with $\mathbb{Z}_m=\mathbb{Z}/m\mathbb{Z}$. 
The main result in \cite{Ham:1998a} is as follows:

\begin{thm} [Hamidoune; cf.~\cite{Ham:1998a}] \label{Hamidoune main}
Let $S$ be a generating subset of $\mathbb{Z}_m$ such that $0 \not \in S$ and $|S| \geq 5$.  Then the number of sums of the subsets of $S$ is at least $\min (m, 2|S| )$.
\end{thm}
As was pointed out in \cite{Ham:1998a}, the result is best possible: if $m=3k$ for some $k \geq 3$, then $$S=\{3,6, \dots, 3(k-1) \} \cup \{1\}$$ has $|\Sigma S|=2k$.  (This example clearly generalizes to noncyclic groups.)

The proof of Theorem \ref{Hamidoune main} in \cite{Ham:1998a} considers two cases: when $2|S| \leq m-1$ and when $2|S| \geq m$.  The proof provided for the first case is correct; in fact, it was delivered for an arbitrary abelian group of order $m$.  However, the author derives the second case from the following:

\begin{lem}  [Hamidoune; cf.~\cite{Ham:1998a}] \label{Hamidoune lemma}
Let $A$ be a subset of $\mathbb{Z}_m \setminus \{0\}$ such that $2|A| \geq m$.  Then $A \cup (2 \hat{\;} A)= \mathbb{Z}_m$.

\end{lem} 
Clearly, if Lemma \ref{Hamidoune lemma} were true, it would immediately yield Theorem \ref{Hamidoune main} in the case when $2|S| \geq m$.  However, Lemma \ref{Hamidoune lemma} is false for every even value of $m$: for example, with $$A=\{1,2, \dots, m/2\},$$ we have $0 \not \in A \cup (2 \hat{\;} A)$.  In fact, when $m \equiv 2$ mod 4, then there are subsets $A$ of $\mathbb{Z}_m$ with the required properties for which $A \cup (2 \hat{\;} A)$ misses two elements of $\mathbb{Z}_m$: for example, for $$A=\{1,2,\dots,(m-2)/4\} \cup \{m/2 , m/2+1, \dots, (3m-2)/4 \},$$ neither $0$ nor $m/2-1$ is in $A \cup (2 \hat{\;} A)$. 

It turns out that for the conclusion of Lemma \ref{Hamidoune lemma}, one must assume that $2|A| \geq m+2$.  More generally, we can prove:

\begin{prop} \label{corrected prop}
Let $G$ be a finite abelian group, and let $G_2$ be the subset---indeed, subgroup---of elements of order at most 2.  
\begin{enumerate}
  \item There is a subset $A \subseteq G \setminus \{0\}$ with $2|A|=|G|+|G_2|-2$ for which $A \cup (2 \hat{\;} A) \neq G$.
\item If $A \subseteq G \setminus \{0\}$ satisfies $2|A| \geq |G|+|G_2|$, then $A \cup (2 \hat{\;} A) = G$.
\end{enumerate}
\end{prop} 
We should point out that $|G|+|G_2|$ is always even, hence our two statements are complementary.  

{\em Proof:}  To prove the first statement, partition $G \setminus G_2$ into disjoint parts $K$ and $-K$ (with $-K$ consisting of the inverses of the elements in $K$).  Then $A=(G_2 \setminus \{0\}) \cup K$ satisfies our requirements.

For our second statement, it suffices to prove that $(G \setminus A) \subseteq 2 \hat{\;} A$.  Let $g \in G\setminus A$ be arbitrary, and let $$L_g=\{x \in G \mid 2x=g\}.$$  We show that if $L_g \neq \emptyset$, then $|L_g| = |G_2|.$
To see this, we choose an element $x \in L_g$, and consider the set $x-L_g$.  (Here and below, for an element $z$  and a subset $Y$ of $G$, we let $z+Y$ denote the set $\{z+y \mid  y \in Y\}$ and $z-Y$ denote the set $\{z-y \mid  y \in Y\}$.) Note that $x-L_g$ has size $|L_g|$ and is a subset of $G_2$, thus $|L_g| \leq |G_2|$.  Similarly, $x+G_2 \subseteq L_g$, so  $|G_2| \leq |L_g|$ as well.

Now let $A_0=A \cup \{0\}$.  Then 
$$|A_0 \cap (g-A_0)|=|A_0|+|g-A_0|-|A_0 \cup (g-A_0)|  \geq 2|A_0| -|G| \geq |G_2|+2.$$  
By the previous paragraph, we then must have an element $a_1 \in A_0 \cap (g-A_0)$ for which $a_1 \not \in L_g$.  Since $a_1 \in g-A_0$, we also have an element $a_2 \in A_0$ for which $a_1=g-a_2$ and thus $g=a_1+a_2$.  But $a_1 \not \in L_g$, and thus $a_2 \neq a_1$.  Now if $a_1=0$, then $a_2 \neq 0$, so $g \in A$, contradicting our assumption.  So $a_1 \in A$ and, similarly, $a_2 \in A$.  Therefore, $g \in 2\hat{\;}A$, as claimed. 
$\Box$

Let us turn now to the proof of Theorem \ref{Hamidoune main}.  We employ the following result:

\begin{thm} [Gallardo, Grekos, et al.; cf.~\cite{GalGre:2002a}] \label{Gallardo, Grekos et al.}
If $m \geq 12$ is even and $|A| \geq  m/2+1$, then $3 \hat{\;} A = \mathbb{Z}_m$.

\end{thm} 

{\em Proof of Theorem \ref{Hamidoune main}:}  As we explained above, we only need to treat the case when $2|S| \geq m$.  In the subcase when $m$ is odd, this inequality is equivalent to $2|S| \geq m+1$; since $G_2=\{0\}$ in this subcase, the second statement of Proposition \ref{corrected prop} implies that $$|\Sigma S| \geq |S \cup (2 \hat{\;} S) |=m=\min\{2|S|,m\}.$$

As the first statement of Proposition \ref{corrected prop} shows, in the subcase when $m$ is even, considering only $S \cup (2 \hat{\;} S)$ is not sufficient.  Luckily, when $m \geq 12$, we can take advantage of Theorem \ref{Gallardo, Grekos et al.}: with $A =S \cup \{0\}$, we have $|A| \geq m/2 +1$, so 
$$|\Sigma S| \geq |3 \hat{\;} S |=m=\min\{2|S|,m\}.$$  This leaves only the cases of $m\in \{6,8,10\}$, which can be checked individually (or see Theorem \ref{Diderrich and Mann} below).  $\Box$

In closing, we mention the following generalization of Theorem \ref{Hamidoune main}:

\begin{thm}  \label{my main}
Let $S$ be a generating subset of an abelian group $G$, and suppose that $0 \not \in S$ and $|S| \geq 5$.  Then the number of sums of the subsets of $S$ is at least $\min (|G|, 2|S| )$.
\end{thm}

Our proof relies on the following 1973 result on the so-called critical number $c(G)$ of $G$ where $$c(G)=\min \{s \in \mathbb{N} \mid A \subseteq G \setminus \{0\}, \; |A| =s \Rightarrow \Sigma A = G \}.$$

\begin{thm} [Diderrich and Mann; cf.~\cite{DidMan:1973a}]  \label{Diderrich and Mann}
Let $G$ be an abelian group of order $2k$ with $k \geq 2$.  
\begin{enumerate}
  \item If $k \geq 5$ or $G \cong \mathbb{Z}_2^3$, then $c(G)=k$.
  \item If $G \cong \mathbb{Z}_4, \mathbb{Z}_6, \mathbb{Z}_8, \mathbb{Z}_2^2,$ or $\mathbb{Z}_2 \times \mathbb{Z}_4$, then $c(G)=k+1$.
\end{enumerate}  
\end{thm}

{\em Proof of Theorem \ref{my main}:}  As we mentioned above, the case when $2|S| \leq |G|-1$ was completed in \cite{Ham:1998a}, so assume that $2|S| \geq |G|$.  If $|G|$ is odd, then, as before, our claim follows from the second statement of Proposition \ref{corrected prop}.  Finally, if $|G|$ is even, the claim follows from Theorem \ref{Diderrich and Mann}.  $\Box$


\begin{thebibliography}{99}


\bibitem{DidMan:1973a} G. T. Diderrich and H. B. Mann, Combinatorial Problems in Finite Abelian Groups.  \emph{A Survey of Combinatorial Theory}, J. N. Srivastava et al., ed., North-Holland (1973).


\bibitem{GalGre:2002a} L. Gallardo, G. Grekos, et al., Restricted addition in $\mathbb{Z}/n\mathbb{Z}$ and an application to the Erd\H{o}s--Ginzburg--Ziv problem.  {\em J. London Math. Soc. (2)} {\bf 65} (2002) 513--523.


\bibitem{Ham:1998a} Y. O. Hamidoune, Adding Distinct Congruence Classes, {\em Combin. Probab. Comput.} {\bf 7} (1998) 81--87. 



\end{thebibliography}
\end{document}